\documentclass{article}
\usepackage{amsmath,amsthm,amsfonts,latexsym}

\begin{document}

\newcommand{\Reali}{{\mathbb{R}}}
\newcommand{\obs}{\mathfrak A}
\newcommand{\fld}{\mathfrak F}
\newcommand{\HH}{\mathfrak H}
\newcommand{\Bi}{\mathfrak B}
\def\C{{\cal C}}
\newcommand{\A}{{\cal A}}
\def\D{{\cal D}}
\def\F{{\cal F}}
\def\H{{\cal H}}
\def\I{{\mathfrak N}}
\def\T{{\cal T}}
\def\R{{\cal R}}
\def\B{{\cal B}}
\def\L{{\cal L}}
\def\P{{\cal P}}
\def\W{{\cal W}}
\def\pL{\psi_\Lambda}

\def\O{{\cal O}}
\def\CC{\mathbb C}
\def\K{{\cal K}}
\def\M{{\cal M}}
\def\N{{\cal N}}
\def\W{{\cal W}}
\def\U{{\cal U}}

\title{
Covariant Sectors with Infinite Dimension \\
and Positivity of the Energy}

\author{
Paolo Bertozzini, Roberto Conti, Roberto Longo\\
{}\\
Dipartimento di Matematica\\
Universit\`a di Roma ``Tor Vergata''\\
Via della Ricerca Scientifica \\
I--00133 Roma, Italy}

\date{April 23, 1997}
\maketitle
\markboth{P. Bertozzini, R. Conti, R. Longo}{Positivity of the Energy}
\renewcommand{\sectionmark}[1]{}
\begin{abstract}
Let $\A$ be a local conformal net of von Neumann algebras on $S^1$ and
$\rho$ a M\"obius covariant
representation of $\A$, possibly with infinite dimension. If $\rho$
has finite index, $\rho$ has automatically positive energy.
If $\rho$ has infinite index, we show the spectrum of the energy always
to contain the positive real line, but, as seen by an example, it may contain
negative values. We then consider nets with Haag duality on $\mathbb
R$, or equivalently
sectors  with non-solitonic extension to the dual net;
we give a criterion for irreducible sectors to have positive energy,
namely this is the case iff there exists
an unbounded M\"obius covariant left inverse.
As a consequence
the class of sectors with positive energy is stable
under composition, conjugation and direct integral decomposition.
\end{abstract}

\vfill
\thanks{E-mail: bertozzi@mat.utovrm.it, conti@mat.utovrm.it,
longo@mat.utovrm.it.\par
Research supported in part by MPI, CNR-GNAFA and INDAM.}

\section{Introduction}

In the past there has been a certain belief that all irreducible
superselection sectors
have finite statistics. Indeed if $\A$ is a translation covariant net on
the Minkowski spacetime and the energy-momentum spectrum has an
isolated mass shell, then, by a theorem of Buchholz and Fredenhagen
\cite{BF}, every positive-energy irreducible representation is localizable in a
space-like cone and have finite statistics.

Nevertheless an analysis of sectors with infinite dimension is possible, in
the context of modular covariant net, as outlined in \cite{GuLoc}, and
this indicated such quantum charges to have a natural occurrence. First
examples of irreducible superselection sectors with infinite dimension have been
constructed, in a simple way, by Fredenhagen \cite{Fr}, associated with
conformal
nets on $S^1$, and moreover there are arguments that in this
context a large natural family of sectors should have infinite dimension
\cite{Re}, thus providing the feeling that infinite statistics might be be
the generic or prevailing situation in low spacetime dimension.

At this point it is natural to begin with a general study of
superselection sectors with infinite dimension. However the extension from the
finite-dimensional to the infinite
dimensional case is certainly far from being straightforward and
requires new methods and insight; it is analogous to the
passage, in the study of group representations, from compact groups to
locally compact groups.

In order to understand the structure of infinite dimensional sectors,
we shall study here the positive energy property. We start with a study
of the finite index case in the one-dimensional conformal
case, with a point of view suitable for generalization.
Classical arguments, see \cite{DHR},
are replaced also due to the failure of Haag duality on the real line
and the occurrence of soliton sectors. Based on modular theory methods,
we show however that the positivity of the energy holds automatically
in the finite index case.

In the context of infinite dimensional sectors we shall then show that
the spectrum of the energy always contains $\mathbb R_+$. But, as we
shall illustrate by a (reducible) example, negative energy values may occur in
general.

We are thus led to characterize the sectors with positive energy.
To this end we study
the basic question whether we may associate
an unbounded left inverse with every
covariant superselection sector with positive energy. We start by
considering a translation-dilation covariant net $\A$ of von Neumann
algebras on $\mathbb R$ obtained by a M\"obius covariant
precosheaf on $S^1$ by cutting the circle at one point. Such
translation-dilation covariant nets are characterized by the
Bisognano-Wichmann geometric action of the modular group of the
von Neumann algebras of half-lines, see \cite{GLW}.

Assuming Haag duality for $\A$ on $\mathbb R$ to hold (for bounded
intervals), we give a positive answer to the above question and in
fact we show that an irreducible sector has positive energy if and
only if it admits a M\"obius covariant unbounded left inverse.

As a consequence we show the class of superselection sectors with
positive energy is closed under composition, conjugation and direct
integral decomposition.

Notice now that the Bisognano-Wichmann dual net $\A^d$ of $\A$
$$
\A^d(a,b)=\A(-\infty,b)\cap\A(a,+\infty)
$$
always satisfies Haag
duality, is conformal thus strongly additive \cite{GLW},
and a covariant representation $\rho$ of $\A$ localized in a bounded
interval $I$ extends to a covariant representation of
$\A^d$, but it may become localized in a half-line, i.e. a soliton
sector \cite{Ro}. However one may construct an extension $\rho_{R}$
of $\rho$ to $\A^d$ localized in a right half-line and an extension $\rho_{L}$
localized in a left half-line. The extensions are easily obtained:
if $I\subset (a,+\infty),$ the restriction of $\rho$ to the C$^*$-algebra
generated by the von Neumann algebras of bounded intervals contained
in $(a,+\infty)$ extends to a normal endomorphism of the its weak
closure $\A(a,+\infty)$, thus, as $a\in\Reali$
is arbitrary, it
gives up an endomorphism of the
C$^*$-algebra $\cup_{a\in\mathbb R}\A(a,+\infty)^{-}$
that restricts to a representation of
the quasi-local C$^*$-algebra of $\A^d$ generated by the local von
Neumann algebras $\A^d(I)$'s, $I$ bounded interval. This representation
is $\rho_R$ and $\rho_L$ is obtained similarly. In general $\rho_R$
and $\rho_L$ are inequivalent representations.

They are equivalent in particular if they are still
localized in a bounded interval, namely the extension of  $\rho$ is
not a soliton.
Our results thus apply to general (non necessarily strongly additive)
M\"obius covariant nets, provided we consider `truly non-solitonic'
covariant sectors.

The idea behind our analysis is to explore the equivalence between the
positivity of the energy and the KMS condition for the dilation
automorphisms of $\A(\Reali_{+})$, that holds in the vacuum sector.
It is rather easy that the latter KMS property passes to charged
sectors, of arbitrary dimension, enabling us to make an analysis by
the Tomita--Takesaki theory.
Our methods rely on an analysis of the unitary representation
of  translation-dilation group, where we give a characterization of
the positive energy representations in terms of domain conditions for
certain associated operators. We then identify these operators with
modular objects furnished by the Tomita-Takesaki theory and use the
domain conditions. In the infinite index case we make use of
Haagerup operator--valued weights, Connes spatial derivatives and Araki
modular operators in particular.
For the convenience of the reader we use \cite{St} as a reference for the
modular
theory.

This paper leaves open the problem whether every irreducible
covariant sector has automatically positive energy. Our work however
indicates that, at least in the strongly additive case, the answer
is likely to be affirmative.

\section{Representations of the translation-dilation group: a criterion for
positive energy}

In this section we analyze the unitary positive-energy  representations
of the dilation-translation group from a point of view
of later use. In particular we are interested in characterizing the
positive energy representations in terms of domain conditions.
To begin with, we relax the hypothesis of positivity
of the energy in the proof of \cite{Lo}, Cor. 2.8.

\medskip
\noindent{\bf 2.1 Lemma.} {\sl
Let $T(a):=e^{iHa},$ $U(t)$
be two $1$--parameter groups on a Hilbert space $\H$
such that $U(t)T(a)U(t)^*=T(e^{-2 \pi t}a)$ for every $t, a \in \Reali.$
Then the spectral projection $P_1$ (resp. $P_2$, $P_3$),
relative to the positive (resp. negative, $0$)
part of the spectrum of $H$
commutes with $T,$ $U$, and thus reduces the representation
on globally invariant subspaces.}\medskip

\noindent{\bf Proof.} {}From the commutation
relation $U(t)T(a)U(t)^*=T(e^{-2 \pi t}a)$,
by differentiation with respect to
$a,$ it follows that $U(t)HU(t)^* = e^{-2\pi t}H.$ For $i=1,2,3$ let $\chi_i$
be the characteristic function of the positive (resp. negative, zero) part
of $\Reali.$ By Borel functional calculus we get:
$\chi_i (U(t) H U(t)^*) = U(t) \chi_i (H) U(t)^* = \chi_i (e^{-2\pi t}H) =
\chi_i (H).$ As $\chi_i (H) = P_i,$ we have the proof.
$\hfill\Box$\medskip

\noindent{\bf 2.2
Lemma.} {\sl Given two $1$-parameter unitary groups
$T(a),$ $U(t)=e^{-2 \pi i t D}$ on the Hilbert space $\H$ such that
$U(t)T(a)U(t)^*=T(e^{-2 \pi t}a)$ for every $t, a \in \Reali,$ then
for every $a \in \Reali$ and $\zeta$ in the
domain of both $e^{-\pi D }$, $e^{-\pi D } T(a),$ we have
$$\| e^{-\pi D } \zeta \|=\|e^{-\pi D } T(a) \zeta \| .$$
The subspace $\D (e^{- \pi D}) \cap \D (e^{- \pi D} T(a))$ of such $\zeta$
is dense for every fixed $a \in \Reali$.}\medskip

\noindent{\bf Proof.}
By symmetry considerations it is sufficient to consider the case $a \geq 0.$
If the generator $-i \frac{d}{da}T(a)|_{a=0}$ of $T$ is non--negative then
$\D (e^{- \pi D}) \subset \D (e^{- \pi D} T(a)),$ $ a \geq 0$
and the thesis is well known
\cite{Lo} (see also Prop. 3.4).
We just outline how to change the arguments
in this reference
in order to have our general statement.

Let us decompose the Hilbert space $\H$ in the direct sum of the two spectral
subspaces $\H_+:=(P_1+P_3)\H ,$ $\H_-:=P_2 \H$.
By Lemma 2.1 these subspaces are invariant for both the
$1$-parameter groups $T(a)$ and $U(t)$, so that it is possible to define
the restrictions $U_+ (t),$ (resp. $U_- (t)$)
and $T_+ (a),$ (resp. $T_- (a)$)
of the $1$-parameter groups to these subspaces.
Now $T_+ (a)$ (resp. $T^{-1}_- (a)$) is a $1$-parameter group with positive
generator satisfying the commutation relation
$U_+ (t)T_+ (a)U_+ (t)^*=T_+ (e^{-2 \pi t}a)$
(resp. $U_- (t)T_-^{-1} (a)U_- (t)^*=T^{-1}_- (e^{-2 \pi t}a)$,)
therefore using
\cite{Lo}, Corollary 2.8, we obtain:
$\| e^{-\pi D_+ } \zeta \|=\|e^{-\pi D_+ } T_+ (a) \zeta \|$
(resp. $\| e^{-\pi D_- } \zeta \|=\|e^{-\pi D_- } T_-^{-1} (a) \zeta \|$)
for every $\zeta$ in the domain of $e^{-\pi D_+ },$
(resp. for every $\zeta$ in the domain of $e^{-\pi D_-}),$
for every $a\geq 0,$ where $D_+$ and $D_-$ are the generators of
the $1$-parameter groups $U_+ (t),$ $U_- (t).$
Let us take $\zeta = (\zeta_+ + \zeta_-)$ in the common domain of $e^{-\pi D}$
and $e^{-\pi D } T(a),$ $a \geq 0,$ then
we have: $\zeta_+$ in the domain of $e^{-\pi D_+ }$ and $\zeta_-$ in the
domain of $e^{-\pi D_-}$ so that we obtain
$\| e^{-\pi D_+ } \zeta_+ \|=\|e^{-\pi D_+ } T_+ (a) \zeta_+ \|$
and
$\| e^{-\pi D_- } \zeta_- \|=\|e^{-\pi D_- } T_-^{-1} (a) \zeta_- \|.$
{}From the second equation (using the fact that
by hypothesis
$T_- (a) \zeta_- $ is in the domain of $e^{-\pi D_-}$), we get
$\| e^{-\pi D_- } T_-(a)\zeta_- \|=\|e^{-\pi D_- } T_-^{-1}(a) T_-(a)\zeta_- \|=
\|e^{-\pi D_- }\zeta_- \|.$ Now, summing the result for the two components
$\zeta_+,$ $\zeta_-,$
from Pitagora's Theorem it is possible to deduce:
$\| e^{-\pi D } \zeta \|=\|e^{-\pi D } T(a) \zeta \|, \ a \geq 0.$
If $a>0$ and $\zeta \in \D (e^{- \pi D}) \subset \D (e^{- \pi D} T(-a)),$
then $T(-a)\zeta \in \D (e^{- \pi D}) \subset \D (e^{- \pi D} T(a)),$
thus
$\| e^{-\pi D } \zeta \| = \|e^{-\pi D } T(a) T(-a)\zeta \| =
\|e^{-\pi D }T(-a)\zeta \|.$
The last part of the statement is now clear.
$\hfill\Box$\medskip

\noindent{\bf 2.3
Corollary.} {\sl Assume that
$T(a),$ and $U(t)=e^{-2 \pi i t D}$ are two $1$--parameter groups as
in Lemma 2.2,
and let $\zeta \in \H$ be a vector such that
$\| e^{-\pi D } \zeta \|=\|e^{-\pi D } T(a) \zeta \| < \infty,$
for some $\ a \geq 0;$ then
$\| e^{-\pi D } \zeta \|=\|e^{-\pi D } T(b) \zeta \|$
for every $0<b<a$
.}\medskip

\noindent{\bf Proof.}
With the same notation as above we know that
$\zeta_+ \in \D (e^{-\pi D_+ })$
and $T_-(a) \zeta_- \in \D (e^{-\pi D_- } ),$
thus
$\zeta_+ \in \D (e^{-\pi D_+ }T(b))$
and
$T_-(b) \zeta_- = T_-(a-b)^* T_-(a)\zeta_-  \in \D (e^{-\pi D_- } ).$
$\hfill\Box$\medskip

We also have the following criterion for the positivity of the energy
\medskip

\noindent{\bf  2.4 Proposition.}
{\sl Given two $1$-parameter groups
$T(a),$ $U(t)=e^{-2 \pi i t D}$ as in Lemma 2.2, the following are equivalent:

\begin{description}

\item{\rm a)} $e^{-\pi D } T(a) \supset T(-a) e^{-\pi D }$, i.e. $T(-a)
e^{-\pi D }$ is hermitean,
for some ($\Leftrightarrow$ for all) $a>0$;

\item{\rm b)} there exists a core $\D_1$ for $e^{-\pi D }$
contained in $\D (e^{-\pi D } T(a)),$
for some ($\Leftrightarrow$ for all) $a>0$;

\item{\rm c)} the generator $H$ of $T(a)$ is positive.
\end{description}

}\medskip

\noindent{\bf Proof.}
The equivalence {\rm a)} $\Leftrightarrow$ {\rm c)} is proved in \cite{Da},
Th. 1.
The implication {\rm c)} $\Rightarrow$ {\rm a)} is also contained in
\cite{Lo}, (proof of) Cor. 2.8 by a different method of proof.
Clearly {\rm a)} $\Rightarrow$ {\rm b)}, thus we have to show that
{\rm b)} $\Rightarrow$ {\rm c)}.
$\D_2 : = e^{-\pi D } \D_1$ is a core for $e^{\pi D },$
and,
for some $a > 0,$
$e^{-\pi D } T(a) e^{\pi D }$ is
isometric on $\D_2$
by Lemma 2.2,
thus by \cite{StZs}, 9.24 the function $t \to U(t)T(a)U(t)^*$
admits an analytic continuation inside the strip
$\{t \in \mathbb C \ | \ - \frac{1}{2} < \Im t < 0 \}$
bounded in norm by $1,$ see also Th. 4.3.
The conclusion may be easily obtained as in
\cite{BGLa}, (proof of) Prop. 2.7;
in fact putting $t=-\frac{i}{4}$ we get $\| T(ia) \| = e^{-a H} \leq 1.$
$\hfill\Box$\medskip

\section{Preliminaries on local conformal precosheaves}

Our analysis  will concern  nets of von Neumann algebras on the
real line. More precisely $\A$ will be a map
$I \mapsto \A(I)$
from the  bounded open intervals $I$ of $\Reali$ to  von Neumann
algebras on a fixed Hilbert space $\H$.
For this net we require the following properties:

\begin{description}

\item{\rm 1)} {\bf Isotony}: $I_1 \subset I_2 \Rightarrow \A(I_1)
\subset \A(I_2)$

\item{\rm 2)} {\bf Locality}: $I_1 \cap I_2 = \emptyset \Rightarrow
[\A(I_1), \A(I_2)] = \{0\};$

\item{\rm 3)} {\bf Covariance}: there exists a strongly continuous
unitary representation $V$ of the translation-dilation group $P$, namely
a semidirect product of $\Reali$ with $\Reali$,
 on $\H$ such that:
 $$V(g)\A(I)V(g)^{-1}= \A(gI), \ g \in P.$$

Here $P$ acts on $\Reali$  ($(a,t)x = a + e^t x$) and we will denote
 by $a, b,\dots \in {\Reali}$  elements of the translation
one--parameter subgroup
 and by $ t,s,\dots\in {\Reali}$  elements of
dilation one-parameter subgroup. We shall frequently  denote
the one-parameter translation (resp. dilation) group simply
by $T(a)$   (resp. $U(t)$).

\item{4)} {\bf Existence of the vacuum}: there exists a unique  (up to
a phase) unit $V$-invariant vector $\Omega\in\H$.

\end{description}

\begin{description}

\item{5)} {\bf Reeh--Schlieder Property}: the vacuum vector $\Omega$ is cyclic
and separating for the von Neumann algebras $\A(I)$.

\item{6)} {\bf Bisognano--Wichmann Property}: the
modular unitary one--parameter group associated (by Reeh--Schlieder Property
and Tomita--Takesaki Theorem) with $(\A({\Reali}_+ ),\Omega)$ coincides
with the rescaled dilation one--parameter
unitary
$$\Delta^{it}_\Omega = U(- 2 \pi t) , \ t \in \Reali.$$
\end{description}
\noindent
Here and in the following, given $S \subset {\Reali},$ we indicate
by $\A_0(S)$ the
$C^*$--algebra generated by all the $\A(I)$,
 with $I \subset S$, and by
$\A(S)=\A_0(S)''$ its weak closure.

In the literature one considers more often M\"obius covariant
precosheaves (also named nets) of von Neumann algebras on the proper
intervals of $S^1$. If one cuts $S^1$ and restricts such a precosheaf
to $\Reali=S^1\backslash \{ \text {point}\}$ one obtains a net on the real
line satisfying the above properties 1 to 6  \cite{BGLa,Bo}. Conversely any
net on
$\Reali$ with the above properties extends uniquely to a M\"obius
covariant precosheaf on $S^1$ \cite{GLW}.

In particular
the modular conjugation $J_{\Reali}$ associated with $\A(\Reali_+),\Omega$
corresponds to the reflection in $\Reali$ with respect to
$0$, thus ``wedge duality'' holds:
$$
\A(a,\infty)'=\A(-\infty,a).
$$
Moreover the generator $H:= -i\frac{d}{da}T(a)|_{a=0}$
of the translation one--parameter group $T$
is  positive \cite{Wi,BGLb}.
In other words positivity of the energy in the vacuum sector is a
consequence of the KMS property for the dilation
automorphism group of $\A({\Reali}_+)$.
We will see how this implication works in
different representations.

Notice that the uniqueness of the vacuum and the positivity of the
energy entail the factoriality of the von Neumann algebras $\A(I)$,
if $I$ is a half-line, thus the irreducibility of the quasilocal
C$^*$-algebra $\A_0(\Reali)$; indeed every net satisfying the above
properties decomposes uniquely into a direct sum of irreducible nets
and irreducibility, uniqueness of the vacuum and factoriality of the von
Neumann algebras of half-lines are equivalent properties.

We shall now consider a morphism $\rho$ of the quasi-local C$^*$-algebra
$\A=\A_0(\Reali)$ localized in a half-line, namely $\rho$ is a
representation of $\A$ on $\H$ such that
$\rho(X)=X$ for every $X \in \A_0(-\infty,a).$
Two such morphisms $\rho, \rho'$ are
said to be equivalent if they are equivalent as representations; thus,
by wedge duality,
there exists a unitary $T \in \A(a,\infty)$ such that
$T\rho(X)=\rho'(X)T$ for every $X \in \A.$

An endomorphism $\rho$ is  covariant if
there exists a unitary
strongly continuous representation $V_\rho : P \to B(\H)$ such that:
$$\rho(\alpha_{g}(X)) = V_\rho(g)\rho(X)V_\rho(g)^{-1},
\ X \in \A, \ g \in P,$$
where $\alpha_g :=\mbox{Ad}(V(g))$.
As far as we consider a covariant irreducible morphism $\rho$ or, more
generally, a finite direct sum of irreducibles (in particular finite
index endomorphisms), the representation
$V_\rho$ providing the covariance is unique,
due to the fact that there are no
non--trivial finite--dimensional representations
of $P.$
In the reducible case, different representations are related
by a cocycle in $\rho(\A)'.$  We shall say that $\rho$ has positive
energy if we can choose $V_{\rho}$ so that the generator of the
translation group is positive.

We shall only consider covariant morphism $\rho$ which are
{\it transportable}, i.e.  localizable in any
half-line $(-\infty,a)$ or $(b,\infty)$. This is in particular the
case of a covariant morphism localizable in an interval.
Thus $\rho$ is normal
on $\A(a,\infty)$ and, if localized in $(a,\infty)$,
extends to a normal endomorphism of $\A(a,\infty)$
denoted $\rho_{(a,\infty)}$ or simply by $\rho$,
if no confusion arises.

By introducing the  notations:
$$\beta^\rho_g := \mbox{Ad}(V_\rho(g)),
\, z_\rho(g):= V_\rho(g)V(g)^*, \quad g \in P,$$
the covariance condition takes the form:

$$\alpha_g \rho \alpha_{g^{-1}} =
\mbox{Ad}(z_\rho(g)^*) \circ \rho \simeq \rho, \ g \in P  \eqno (3.1)
$$

\noindent
and $z_\rho(g)$
satisfies the following $\alpha$--cocycle identity:
$$z_\rho(gg')=z_\rho(g)\alpha_g(z_\rho(g')),
\quad g,g' \in P .
$$
Let $\rho$ be localized in the half-line $I=(a,\infty)$ then
$\alpha_g \rho \alpha_{g^{-1}}$ is
localized in $gI$ and from formula (3.1) we obtain $z_\rho(g) \in
\A((\tilde I)')'=\A(\tilde I)$ where
$\tilde{I}$ is the largest half-line between $I$ and $gI$.

We will often use the notations:
$$M:=\A(0, +\infty), \ M_b:=\alpha_b(M)=\A(b,+\infty), \
M_b^\rho:=\beta^\rho_b(M), \ b \in {\Reali}_+.$$
As $z_\rho(b) \in M_{a}$, $b>0$ it follows that
$$ M_{b}=M^\rho_{b},\quad b<a.
$$
The Bisognano--Wichmann Property states that the one--parameter group
$t \mapsto \alpha_{-2\pi t}, \ t \in {\Reali},$ coincides on $M$ with the
modular group $t \mapsto \sigma_t:=\mbox{Ad}(\Delta^{it}_\Omega), \ t \in
{\Reali}$ and, since the cocycle $ z_\rho(- 2 \pi t)$ is localized in $M,$ by
Connes' Theorem there exists a unique semifinite normal faithful (s.n.f.)
weight $\psi_\rho$ on $M$ whose Radon--Nikodym derivative with respect
to the vacuum state $\omega:= (\Omega, \cdot \Omega)$ is given by
$(D\psi_\rho : D\omega)_t=z_\rho(-2 \pi t),$ see \cite{St}, Sect. 11.
Then $t \mapsto \beta^\rho_{-2\pi t}=\mbox{Ad}(z_\rho(-2\pi
t))\circ \alpha_{-2\pi t}, \
t \in {\Reali}$, is the modular group associated to the weight $\psi_\rho$
on $M.$

\section{Automatic positivity of the energy in the finite index case}

Although we shall be mainly interested in sectors with infinite
dimension, our proof will be more transparent by a previous analysis
of the finite index case
(recall that the index is the square of the dimension, see \cite{L}).

Let $\rho$ be an endomorphism of $\A$ localized in $I \subset {\Reali}_+ .$
In the finite index case the following analog of the Kac-Wakimoto formula holds
\cite{Lo}:
$$(D \varphi_\rho : D \omega)_t =
d(\rho)^{-it} z_\rho(-2\pi t),$$
where  $\omega =
(\Omega,\cdot\Omega)$ is the restriction of the vacuum state
to $M=\A({\Reali}_+)$, $\varphi_\rho$ is the state
$\omega \circ \phi_\rho$
on  $M$ and
$\phi_\rho = \rho^{-1} \circ E_\rho $ is
the minimal left inverse of $\rho: M \to M,$
with $E_\rho: M \to \rho(M)$ is the minimal
conditional expectation.

We recall the following:

\medskip
\noindent{\bf 4.1 Proposition. \cite{Lo}} {\sl
Let us assume that $\rho$ is
covariant
with finite index as above,
and let $\psi_\rho$ be a positive linear functional on $M;$
then the following are equivalent:
\begin{description}

\item{\rm a)} $\psi_\rho$ is normal, faithful and
$(D\psi_\rho: D\omega)_t = z_\rho(- 2 \pi t),$ $t \in \Reali$

\item{\rm b)} $\psi_\rho = d(\rho) \ \omega \circ \phi_\rho,$

\item{\rm c)} $\psi_\rho(XY^*)=
(e^{-\pi K_\rho}X\Omega,e^{-\pi K_\rho}Y\Omega),$ $X, Y \in M,$

\item{\rm d)}
$\psi_\rho$ is normal faithful,
$\sigma^{\psi_\rho}_t \circ \rho = \rho \circ \sigma^{\omega}_t,$
$t \in \Reali,$ and $\psi_\rho | \rho(M)'\cap M$ is a trace whose value on a
central projection $p$
is $\psi_\rho(p) = d(\rho_p),$ where $\rho_p$ is the subrepresentation
associated
to $p.$
(In particular, if $\rho$ is irreducible, this last condition
reduces to $\psi_\rho(I)=d(\rho)$).
\end{description}
}\medskip

\noindent{\bf Proof.}
We sketch the first part of the proof.
We assume $\rho$ to be irreducible.\par
We consider the states $\omega,$ $\varphi_\rho:=\omega \circ \phi_\rho$
on the von Neumann algebra $M.$
Note that $\varphi_\rho \circ \rho = \omega$ i.e.
$\varphi_\rho |_{\rho(M)} = \omega \circ \rho^{-1}.$
{}From $\varphi_\rho = \varphi_\rho \circ E_\rho$
it follows that
$\rho(M)$ is $\sigma^{\varphi_\rho}$--stable by Takesaki's theorem,
therefore
$$\sigma^{\varphi_\rho}_t |_{\rho(M)} =
\sigma^{{\varphi_\rho}|_{\rho(M)}}_t =
\rho \circ \sigma^{\omega}_t \circ \rho^{-1} .$$
Now defining $v_t := (D \varphi_\rho : D \omega)_t \in M$ then we have
$v_t \sigma^{\omega}_t (X) v^*_t = \sigma^{\varphi_\rho}_t(X), X \in M,$
thus
$v_t \sigma^{\omega}_t (\rho(X)) v^*_t = \sigma^{\varphi_\rho}_t(\rho(X)) =
\rho \circ \sigma^{\omega}_t \circ \rho^{-1} \rho(X) =
\rho \sigma^{\omega}_t(X) = \beta^\rho_{-2 \pi t} \rho(X) =
z_\rho(-2 \pi t) \sigma^{\omega}_t (\rho(X)) z_\rho(-2 \pi t)^*, X \in M.$
Hence
$$z_\rho(-2 \pi t)^* v_t \in \sigma^{\omega}_t (\rho(M))' \cap M =
\sigma^{\omega}_t (\rho(M)' \cap M) = \mathbb C.$$
Now to complete the argument as regard to the phase $d(\rho)$ in b), we
refer to \cite{Lo} part 1.
The proof of the point c) will be easily obtained by polarization of
the first formula contained in Proposition 4.2.

Assuming d), to obtain a), notice that this condition determines
$(D\psi_\rho : D\omega)_t$ up to the multiplication by a cocycle
in $\rho(M)'\cap M$ hence $\psi_\rho$ is determined by the
specification of $\psi_\rho | \rho(M)'\cap M$ that we require to be
$d(\rho)$--times the restriction to $\rho(M)'\cap M$ of the minimal
expectation of $M$ onto $\rho(M).$
$\hfill\Box$\medskip

Now we have
$\beta^\rho_{-2 \pi t} = \sigma^{\psi_\rho}_t (= \sigma^{\varphi_\rho})
=\mbox{Ad}(\Delta^{it}_\xi)=\mbox{Ad}(\Delta^{it}_{\xi,\Omega})$ on $M$,
where $\varphi_\rho=(\xi, \cdot \xi),$ $\| \xi \|=1$ and $\xi$ is cyclic
for $M$
(e.g. $\xi$ is the vector representative of $\varphi_\rho$
in the natural cone of $M$ given by $\Omega$).
Clearly $\psi_\rho \circ \beta^\rho_t = \psi_\rho.$

Using the fact that $\rho$ is localized in ${\Reali}_+$
and recalling the definition of $\psi_\rho$
it is easy to check that
$U_\rho(-2 \pi t) = z_\rho(-2\pi t)U(-2\pi t)=
z_\rho(-2\pi t)\Delta^{it}_{\Omega}$
coincides up to the phase $d(\rho)^{it}$ with
$\Delta^{it}_{\xi,\Omega}=\Delta^{it}_{\varphi_\rho,\omega},$
where $\Delta^{it}_{\xi,\Omega}$ is the Araki relative modular
operator, see \cite{BrRo}, namely
$S_{\xi,\Omega}=J_{\xi,\Omega}\Delta^{\frac{1}{2}}_{\xi,\Omega}$
is the polar decomposition of the closure of
$X \Omega \to X^* \xi, X \in M$.
In fact we have
$$
2 \pi K_\rho
= - \mbox{log} \Delta_{\xi,\Omega} - \mbox{log}d(\rho),
$$
see \cite{Lo}.
Hence by the commutation relations of the group $P$ we obtain
$$
\Delta^{it}_{\xi,\Omega} T_\rho(a) \Delta^{-it}_{\xi,\Omega} =
T_\rho(e^{- 2 \pi t}a), \ t, a \in \Reali .
$$
The argument given
in Lemma 2.2 amounts to a proof that the invariance
condition
holds if the dimension of $\rho$ is finite, cf. \cite{Lo}, Prop. 2.11.

\medskip

\noindent{\bf 4.2 Proposition.} {\sl Let $\rho$ be a covariant endomorphism
with finite dimension localized in $I \subset {\Reali}_+,$
and let $M,$ $U_\rho(a), \ a \in \Reali,$ $\varphi_\rho$ be as above.
Then $\varphi_\rho$ is $\beta^\rho_a$--invariant
on $M,$
for every $a \geq 0$ .}
\medskip

\noindent{\bf Proof.}
If $X \in M$ we have
\begin{align*}
\varphi_\rho(X^*X)&=(\xi,X^*X\xi)=\|X\xi\|^2=\|S_{\xi,\Omega}X^*\Omega\|^2\\
&=\| J_{\xi,\Omega}\Delta^{\frac{1}{2}}_{\xi,\Omega} X^* \Omega\|^2=
\| \Delta^{\frac{1}{2}}_{\xi,\Omega} X^* \Omega\|^2
\end{align*}
(see formula (3.3) in \cite{Lo}), therefore if $a >0$
and $z_{\rho}(a)$ denotes
the cocycle with respect to the translation by $a$,
\begin{align*}
\varphi_\rho (T_\rho(a)X^*XT_{\rho}(a)^*)
&=\| \Delta^{\frac{1}{2}}_{\xi,\Omega} T_\rho(a) X^* T_{\rho}(a)^* \Omega\|^2\\
&=\| \Delta^{\frac{1}{2}}_{\xi,\Omega} T_\rho(a) X^* z_\rho(-a)T(a)^* \Omega\|^2
\\
&=\| \Delta^{\frac{1}{2}}_{\xi,\Omega} T_\rho(a) X^* z_\rho(-a) \Omega\|^2\\
&=\| T_\rho(a)\Delta^{\frac{1}{2}}_{\xi,\Omega}X^* z_\rho(-a) \Omega\|^2
\ \ \mbox{(by Lemma 2.2)}
\\
&=\| \Delta^{\frac{1}{2}}_{\xi,\Omega}(X^* z_\rho(-a)) \Omega\|^2
=\| z_\rho(-a)^* X \xi \|^2 = (X\xi,X\xi).
\end{align*}
$\hfill\Box$

\medskip

In the sequel we write $T \ \eta \ M$ to denote that the (generally unbounded)
linear operator $T$ on $\H$ is affiliated to $M \subset B(\H)$,
see \cite{StZs}, 9.7.

\medskip
\noindent{\bf 4.3 Proposition.}
{\sl
Let $\rho$ be as above.
For any given $a \geq 0$,
the function
$t \mapsto \Delta^{it}_{\xi,\Omega} T_\rho(a) \Delta^{-it}_{\xi,\Omega} =
T_\rho(e^{-2 \pi t}a)$ admits an analytic continuation inside the strip
$\{z \in \mathbb C \ | - \frac{1}{2} < \Im z <0 \}$
which is bounded in norm by $1.$
}\medskip

\noindent{\bf Proof.}
The existence of the analytic continuation inside the strip
follows from general arguments, see e.g. \cite{BDL} p. 241.
The bound $1$ is thus a consequence of Hadamard three line theorem,
once we check
that the norm of the function is bounded by $1$ on the lines $z=0$ (this
is obvious),
$z=-\frac{i}{2} + t, t \in \Reali$ and has a priori global bound on the
entire strip.

We now check that
$\Delta^{\frac{1}{2}}_{\xi,\Omega} T_\rho(a)
\Delta^{- \frac{1}{2}}_{\xi,\Omega}$, or eqivalently
$S_{\xi,\Omega}T_{\rho}(a)S_{\xi,\Omega}^{-1}$,
is extended by an
isometric operator.
We write for short
$S=S_{\xi,\Omega}.$\par

Note that $S \ T_\rho(a) \ S^{-1}$ is isometric
on the dense subspace $M \xi:$
given $X \in M,$
we have
\begin{align*}
S \ T_\rho(a) \ S^{-1} \ X^* \xi & = S T_\rho(a) X \Omega  =
S \ z_\rho(a) T(a) X T(a)^* \Omega \\
& =  T(a) X^* T(a)^* z_\rho(a)^* \xi = T(a) X^* T_\rho(a)^* \xi
\end{align*}
therefore $M \xi \subset \D (S T_\rho(a) S^{-1}) \subset \D (S^{-1})$ and
$$
\| \ S \ T_\rho(a) \ S^{-1} \ X^* \xi \| = \| X^* \xi \|, \ X \in M
$$
by the $\beta^\rho_a$--invariance of $\varphi$ showed in Prop. 4.2.

As is known
$\D (S) =\{T \Omega \ | \ T \ \eta \ M, \ \Omega \in \D (T), \ \xi \in \D
(T^*)\},$
and using
Lemma 2.2
it is direct to verify that
$\| \ S \ T_\rho(a) \ S^{-1} \ \zeta \| = \| \zeta \|$ for every
$\zeta \in \D (S^{-1})$ for which the l.h.s. is well defined.

It remains to check the a priori bound on the strip. This may be
derived by the bound on matrix
coefficients $|(\Delta^{iz}_{\xi,\Omega} T_\rho(a)
\Delta^{-iz}_{\xi,\Omega}\zeta_1,\zeta_2)|\leq
\|\zeta_1\|\|\zeta_2\|$, with $ - \frac{1}{2} < \Im z <0 $, and
$\zeta_1,\zeta_2$
in spectral subspaces for $\log\Delta_{\xi,\Omega}$ with respect to
bounded intervals. The bound follows by the three line
theorem, because the involved fuctions are bounded on the strip.
$\hfill\Box$

\medskip
Notice that the same result may be obtained using
the Prop. 2.4 or  the Proposition stated in \cite{StZs}, p. 219.

\medskip
\noindent {\bf 4.4 Corollary.}
{\sl
Let $\rho$ be a covariant endomorphism with finite dimension localized
in $I \subset {\Reali}_+,$
and $M,$ $T_\rho(a)$ as above. Then in the sector
$\rho$ the energy (the generator of the $1$--parameter group $T_\rho(a)$)
is positive.}

\medskip
\noindent {\bf Proof.}
It is immediate from (the proofs of) Prop. 2.4
(the core is $M \Omega$)
and Prop. 4.2, cf. Prop. 4.3.
$\hfill\Box$

\section{Infinite index (weight) case }
\subsection{General considerations}
We shall now begin an
analysis of sectors with infinite dimension, by an extension of the
previous methods.

In the following $\A$ is a net on $\Reali$ as in the previous section,
namely $\A$ is obtained by restricting a local conformal
precosheaf of von Neumann algebras on $S^1$.

Let $\rho$ be a covariant endomorphism localized in a half-line
strictly contained in $\Reali_+$ as before, but
not necessarily with finite index.  We know from \cite{Lo}, Sect. 2, that
$$
U_\rho(- 2 \pi t) =: e^{- 2\pi i t K_\rho} =
\Delta(\psi_\rho /\omega')^{it}   \eqno (5.1)
$$
(cf. the paragraphs following Prop. 4.1 and Lemma A.2)
where $\omega'$ is the restriction of the vacuum state
to $M',$ and $\Delta(\psi_\rho /\omega')$ is the Connes' spatial derivative
with $\psi_\rho$ the weight on $M$ defined at the end of section 3.

Although $\psi_\rho$ is unbounded in general,
$\omega'$ is a state, represented by the
vector $\Omega.$
This suggests  the formula
$$
\psi_\rho(XX^*)=\|e^{- \pi K_\rho} X\Omega\|^2
$$
still to express $\psi_\rho$
and  to be useful in proving its invariance properties,
cf. \cite{Lo} (3.3).
We shall show this is in fact the case. To this end
we need to discuss certain aspects
of the spatial theory of von
Neumann algebras, that
may be of independent interest, and we collect them
in  Appendix \ref{Appendix A}. We refer to this appendix for
notations and notions used here below.

Notice now that equation $(5.1)$ gives the dilation--translation
commutation relations
$$ \Delta(\psi_\rho /\omega')^{it} T_\rho(a) \Delta(\psi_\rho
/\omega')^{-it}
= T_\rho(e^{- 2 \pi t}a), \ a, t \in \Reali,$$
so we may still apply the analysis made in Section 2.

We use the following notations:
$$\psi_a:=\psi_\rho \circ \mbox{Ad}(T_\rho(a)),\quad a>0,$$
is a faithful normal weight on $M;$
$$\I:=\I_{\psi_\rho},\quad \I_a:=\{X \in M \ \  | \ \ \psi_a(X^*X) < \infty
\},$$
are left ideals in $M$ for any $a > 0.$

Note that $\sigma_t^{\psi_\rho}(\I_a)=\I_{e^{2\pi t}a},$ $a>0,$ $t \in \Reali$
and that if $X \in M,$ and  $a>0$ then, cf. Prop. 4.2 and Lemma A.1,
$$X^* \Omega \in \D (\Delta(\psi_\rho /\omega')^{\frac{1}{2}}
T_\rho(a))
\Leftrightarrow
X \in \I_a .$$
As a first result we show that the spectrum of
the generator $H_\rho$ of the translation
group relative to
the sector $\rho$ always contains the positive real line.
\medskip

\noindent{\bf 5.1 Proposition.}
{\sl
Let $\A$ be a local net of von Neumann algebras on $\Reali$ as in Sect. 3
and $\rho$ a covariant morphism localized in $(1,+\infty).$
Then $[0,+\infty) \subset \mbox{Sp}(H_\rho).$}\medskip

\noindent{\bf Proof.}
We already know that $\I$ is dense in $M$
by the semifiniteness of $\psi_\rho.$
By dilation covariance $\mbox{Sp}(H_\rho)$ is either $\Reali_+ \cup \{0\},$
$\Reali_- \cup \{0\},$ or $\Reali$
thus we assume it to be $\Reali_-\cup \{0\}$ to find a contradiction.
We then have
$\D (\Delta(\psi_\rho /\omega')^{\frac{1}{2}}) \subset
\D (\Delta(\psi_\rho /\omega')^{\frac{1}{2}} T_\rho(a)^*),$ $a >0,$
see Sect. 2.
We take $X \in \I^*$ and define $\zeta = T_\rho(a)^* X T(a) \Omega,$
thus we have $\zeta \in \D (\Delta(\psi_\rho /\omega')^{\frac{1}{2}}).$
Now we observe that the closed operator $T_\zeta \ \eta \ M$ as defined in
Lemma A.3
coincides with the bounded operator $T_\rho(a)^* X T(a)$
which a priori is in $ M_{-a};$
in fact the two operators coincide when restricted to the dense
vector space $M_{-a}' \Omega.$ Therefore $T_\zeta$ is bounded and
$T_\zeta = T_\rho(a)^* X T(a) \in M.$
It follows that
$X \in M_a$
whenever $a >0$ is small enough.
Therefore we have $\I \subset M_a$ thus $M \subset M_a$
which is not possible.
$\hfill\Box$\medskip

\noindent{\bf 5.2 Lemma.}
{\sl $0 < b < a$ $\Rightarrow$ $\I \cap \I_a \subseteq \I \cap \I_b.$
}\medskip

\noindent{\bf Proof.}
By Lemma A.1 if $X^*\in  \I \cap \I_a,$ then
$X \Omega \in \D (\Delta(\psi/ \omega)^{\frac{1}{2}}) \cap $
$\D (\Delta(\psi/ \omega)^{\frac{1}{2}} T_{\rho}(a)),$ thus by Corollary 2.3,
$X \Omega \in \D (\Delta(\psi/ \omega)^{\frac{1}{2}}) \cap
\D (\Delta(\psi/ \omega)^{\frac{1}{2}} T_{\rho}(b)),$ hence $X^* \in  \I
\cap \I_b$ by Lemma A.1.
$\hfill\Box$\medskip

\noindent{\bf 5.3 Proposition.}
{\sl
$\I \subset \I_a$ for any $a>0$ if and only if $\I = \cup_{a>0} \I \cap \I_a.$
}\medskip

\noindent{\bf Proof.}
One implication is obvious. To show the converse assume $\I = \cup_{a>0} \I
\cap \I_a$
and let $X \in \I$ be fixed. Then there exists $a=a_1 >0$ such that $X \in
\I \cap \I_a$
i.e. $\psi_\rho \circ (X^*X) = \psi_\rho \circ \beta^\rho_a(X^*X)< \infty$, thus
$\beta^\rho_a(X) \in \I;$
by iteration we find an increasing sequence $a_n >0$ such that
$\psi_\rho \circ \beta^\rho_{a_n} (X^*X) < \infty$.

Let us consider $a_*:=\mbox{sup}\{a \in \Reali_+ \ | \ X \in \I_a\}$,
$a_n \nearrow a_*.$
But $\psi_\rho \circ \beta^\rho_{a_*}(X^*X) \leq
\underline{\mbox{lim}}_n \psi_\rho \circ \beta^\rho_{a_n}(X^*X) =
\psi_\rho(X^*X)$
by lower semicontinuity,
therefore $a_*=\infty$ and we are done.
$\hfill\Box$\medskip

\noindent{\bf 5.4 Proposition.}
{\sl We have
$\cap_{a>0} \I_a \subseteq \I.$
}\medskip

\noindent{\bf Proof.}
If $X \in \cap_{a>0} \I_a$,
we have
$\psi_\rho \circ \beta^\rho_a(X^*X) = \psi_\rho \circ \beta^\rho_b(X^*X)$
for every $a, b >0.$
$\hfill\Box$\medskip

We summarize the last results in the following proposition, although
only the equivalence between the first two points will be needed.
\medskip

\noindent{\bf 5.5 Proposition.}
{\sl The following assertions are equivalent:
\begin{description}
\item{1)} $\rho$ has positive energy;
\item{2)} $\I \subset \I_a$ for some ($\Leftrightarrow$ for all) $a > 0;$
\item{3)} $\I = \cup_{a > 0}(\I \cap \I_a);$
\item{4)} $\I = \cap_{a > 0} \I_a.$
\end{description}
}

\noindent{\bf Proof.} We shall show the equivalence $1) \Leftrightarrow
2),$
the other results are immediate from Lemmata 5.2, 5.3, 5.4.

\noindent
$1)\Rightarrow 2)$ {}From Lemma A.1 we have:
$\I^* = \{ X \in M \ | \ X\Omega \in
\D (\Delta(\psi / \omega ')^{\frac{1}{2}})\}$ and
$\I_a^* = \{ X \in M \ | \ X\Omega \in
\D (\Delta(\psi / \omega ')^{\frac{1}{2}} T_\rho (a)) \},$ cf. Prop. 4.2.
Now, as already recalled in the proof of
Lemma 2.2 (from \cite{Lo}, Corollary 2.8), if $\rho$ has positive energy,
$\D (\Delta(\psi / \omega ')^{\frac{1}{2}}) \subset
\D (\Delta(\psi / \omega ')^{\frac{1}{2}} T_\rho (a))$ for any $a > 0$ and
so: $\I^* \Omega \subset \I_a^* \Omega.$

\noindent
$2) \Rightarrow 1)$ {}From \cite{St}, p. 94--95, we have that $\I^* \Omega$
is a core for $\Delta(\psi / \omega ')^{\frac{1}{2}}.$
If $\I^* \subset \I_a^*,$ by Lemma A.1
$\I^* \Omega \subset \D (\Delta(\psi / \omega ')^{\frac{1}{2}} T_\rho (a))$
and so by Proposition  2.4 b), c), we get the positivity of the energy.
$\hfill\Box$

\subsection{Sectors of Haag dual nets on $\Reali$}
In the following $\A$ is a net of von Neumann algebras on $\Reali$
satisfying the six properties listed in Sect. 3.
Furthermore
we require Haag duality on $\Reali.$
Equivalently we assume our net to be strongly additive, meaning that
$\A(a,b) \vee \A(b,c) = \A(a,c)$ for every $a < b <c,$
$a,b,c \in \Reali.$
It follows that $\A(a,b) ' \cap \A(a,c) = \A(b,c)$
\cite{GLW}.
Let $\rho$ be  a covariant morphism of $\A$  localized in
$(b,+\infty),$ $0 < b,$. As always we assume that $\rho$ is
localizable in every half-line, thus it extends to a normal
endomorphism of $\A(b,\infty)$. As recalled in Section 3, we obtain a localized
cocycle $t \to z_\rho(- 2 \pi t) \in M:=\A(0,\infty)$, thus
a s.n.f. weight $\psi_\rho$ on $M$ such that
$(D \psi_\rho: D \omega)_t = z_\rho(- 2 \pi t), t \in \Reali.$
We have $U_\rho(-2 \pi t)=\Delta(\psi_\rho /\omega')^{it}$
(see formula (5.1)).

For each $X \in M$ we have
\begin{align*}
\sigma^{\psi_\rho}_t (\rho(X)) & =
z_\rho(- 2 \pi t) \alpha_{- 2 \pi t}(\rho(X)) z^*_\rho(- 2 \pi t)\\
&= z_\rho(- 2 \pi t) ( \alpha_{- 2 \pi t} \rho \alpha_{2 \pi t}
(\alpha_{- 2 \pi t}(X)) ) z^*_\rho(- 2 \pi t)
= \rho(\alpha_{- 2 \pi t}(X)) \\
&= \rho \alpha_{- 2 \pi t} \rho^{-1} (\rho(X)) =
\sigma^{\omega \rho^{-1}}_t (\rho(X)), \ t \in \Reali,
\end{align*}
therefore by Haagerup's theorem
\cite{St}, p. 164,
there exists a unique s.n.f.
operator valued weight
$E:=E_\rho: M^+ \to \overline{\rho(M)}^+$
(where the symbol $\overline{\rho(M)}^+$ denotes the extended positive part
of the
von Neumann algebra $\rho(M)$,
see Appendix \ref{Appendix A})
such that $\psi_\rho = \omega \rho^{-1} E.$
Here $\omega \rho^{-1}$ has to be thought of as the unique extension
to $\overline{\rho(M)}^+$ such that
$\omega \rho^{-1}(n) = n(\omega \rho^{-1}),$ $n \in \overline{\rho(M)}^+.$
In particular $E(X^* Y X) = X^* E(Y) X,$ for every $X \in \rho(M), Y \in M^+.$
$E$ may be uniquely ``extended'' to a
(densely defined)
linear mapping (also denoted by $E$)
${\cal M}_E := \mbox{lin} \{ X \in M^+ ; E(X) \in \rho(M)^+ \} \to \rho(M)$
with image a ultraweakly dense two--sided ideal in $\rho(M)$,
such that
$E(\rho(Y)X\rho(Z))=\rho(Y)E(X)\rho(Z) \in \rho(M)$ for every $Y, Z \in M,$
$X \in {\cal M}_E$.
Clearly ${\cal M}_E \subseteq {\cal M}_{\psi_\rho}.$
We consider the unbounded left inverse of $\rho$ defined by
$\phi_\rho:=\rho^{-1} \circ E : M^+ \to \rho^{-1}(\overline{\rho(M)}^+)
= \overline{M}^+;$
by linearity $\phi_\rho: {\cal M}_E \to M.$
Clearly $\rho \phi_\rho = E,$ $\psi_\rho = \omega \phi_\rho$ both on $M^+,$
${\cal M}_E,$ and
$\phi_\rho(\rho(Y)X\rho(Z)) = Y \phi_\rho(X) Z$ for every $Y, Z \in M,$
$X \in {\cal M}_E$.
\medskip

\noindent{\bf 5.6 Lemma.} {\sl
With the notations above we have $E(X) \in \overline{M_a}^+$
for every $X \in {M_a}^+$ with $0<a$ sufficiently small.
}\medskip

\noindent{\bf Proof.}
Let $N:=M_a'\cap M,$ making use of the strong additivity property
and duality we obtain that $N=\A(0,a)$ (it is immediate to see that
$\A(0,a) \subset N,$ on the other side if $x \in N$ then $x \in
\A(0,a)=\A(0,a)''=(M' \vee M_a)'$).
If $a$ is sufficiently small, so that
$\rho$ is localized in $(a,\infty)$
then $\rho(u)=u$ for every $u \in N.$
For every $X \in M_a^+$ and unitary $u \in N$
we have $X=uXu^*,$ so that $E(X)=E(uXu^*).$
Now using the fact that $u=\rho(u)$ (because $u \in N$), and that $E$
is an operator valued weight, it follows that $E(X)=uE(X)u^*$
(in fact $E(uXu^*)=E(\rho(u)X\rho(u)^*)=\rho(u)E(X)\rho(u)^*=uE(X)u^*$).
We know that $E(X)$ can be uniquely written as $he \ + \infty (I - e),$
with $e \in P(\rho(M)) \subset P(M)$
($P(M)$ is the set of all the projections of $M$)
and $h \ \eta \ e\rho(M)e$ positive.
On the other side $E(X)=uE(X)u^*=uhu^* ueu^* \ + \infty (I - ueu^*)$ and
using the uniqueness of the decomposition
we obtain: $e=ueu^*$ and $h=uhu^*$ for every $u \in N.$
{}From this, with the help of strong additivity ($N' \cap M=M_a$),
we obtain $e \in M \cap N'=M_a$ (because $e \in P(M)\subset M$ and
$e$ commutes with every $u \in N$). In the same way we can say that
$h \ \eta \ e{M_a}e:$
the bounded parts of $h$ are in $e{M_a}e,$
but the bounded parts $h_n$ of $h$ are in $N'$ (because $h$ commutes
with every unitary in $N$) and in $eMe$ (because $h \ \eta \ eMe$); using
the fact that $eMe \cap N'= e{M_a}e$ (clearly $e{M_a}e \subset eMe \cap N'$
since $e \in M_a \subset N',$ and if $x \in eMe \cap N',$ $x=eme=e(eme)e$ with
$eme \in M \cap N'= M_a$) we get $h_n \in e{M_a}e$ and
thus $h \ \eta \ e{M_a}e$. {}From the fact that $E(X)= he \ + \infty (I - e)$
with $h \ \eta \ e{M_a}e$ and $e \in M_a$ it follows the result:
$E(X) \in \overline{M_a}^+.$
$\hfill\Box$\medskip

\noindent{\bf 5.7 Lemma.} {\sl
Let $M,$ $N$ be von Neumann algebras on the Hilbert space $\H.$
We have $\overline{M}^+ \cap \overline{N}^+ = \overline{M \cap N}^+.$
}\medskip

\noindent{\bf Proof.}
Both $\overline{M}^+,$ $\overline{N}^+$ can be embedded in $\overline{B(\H)}^+.$
$m \in \overline{M}^+ \cap \overline{N}^+$ can be uniquely written as
$h e \ + \infty (I -e), $ with $e \in P(M),$ $h \ \eta \ eMe,$
and
$h' e' \ + \infty (I -e'), $ with $e' \in P(N),$ $h' \ \eta \ e'Ne'.$
It follows that $e=e' \in P(M \cap N),$ and
$h=h' \ \eta \ M \cap N.$
$\hfill\Box$
\medskip

\noindent{\bf 5.8 Lemma.} {\sl
Let $M = \A(0,\infty),$ and $\rho$  localized in $(b,+\infty),$ $0<b.$
Then we have
$$\beta^\rho_a(M) \cap \rho(M)=\rho(\beta^\rho_a(M)) $$
whenever $0<a<b.$
}\medskip

\noindent{\bf Proof.}
The inclusion
$\beta^\rho_a(M) \cap \rho(M) \supset \rho(\beta^\rho_a(M))$ is obvious.
We choose $u_n \in (\rho,\rho_n)$ unitaries with $\rho_n$ localized in
$(-\infty,c_n)$ with $c_n \to -\infty.$
Any weak limit point of the sequence $\mbox{Ad}(u_n)$ is a map
$\tilde{\phi} : B(\H) \to B(\H)$
such that $\tilde{\phi} \rho = \mbox{id}$ on $\cup_{l \in \Reali}\A(l,\infty).$
$\tilde{\phi}$ is a non normal left--inverse of $\rho,$
cf. \cite{DHR}.
Let $X \in \beta^\rho_a(M) \cap \rho(M).$
Then $X = \rho(Y)$ for some $Y \in M,$ and $Y = \tilde{\phi}(X).$
For every $Z \in \A(I),$ $I \subset (-\infty,a)$ we have
$ZY=
Z\tilde{\phi}(X)=\tilde{\phi}(\rho(Z)X)=
\tilde{\phi}(ZX)=\tilde{\phi}(XZ)=YZ$
therefore $Y$ commutes with $\cup_{I \subset (-\infty,a)} \A(I).$
By duality $Y \in \beta^\rho_a(M)=\alpha_a(M)$ and we are done.
$\hfill\Box$\medskip

\noindent{\bf 5.9 Proposition.} {\sl
Let $\rho$ be a covariant (transportable) morphism of $\A$
localized in $(b,+\infty),$ $b>0,$
$M := \A(0,\infty),$
$E: M^+ \to \overline{\rho(M)}^+$ defined as above, and
$E_a:= \rho \alpha_{-a} \rho^{-1} E \beta^\rho_a.$
Then $E_a: M^+ \to \overline{\rho(M)}^+$ is a s.n.f. operator valued weight,
for $0<a<b$.
If $\rho$ is irreducible and $\psi_a$ is
semifinite, we have $E =E_a,$ thus $\psi = \psi_a.$
}\medskip

\noindent{\bf Proof.}
$X \in M^+$ $\Rightarrow$ $\beta^\rho_a(X) \in {M_a}^+$ $\Rightarrow$
$$\begin{array}{ccc}
E\beta^\rho_a(X) \in \overline{\rho(M)}^+ \cap \overline{M_a}^+ & &
\mbox{(by Lemma 5.6)} \\
= \overline{\rho(M) \cap M_a}^+  & &
\mbox{(by Lemma 5.7)} \\
= \overline{\rho(M_a)}^+ & &
\mbox{(by Lemma 5.8)}
\end{array}$$
$\Rightarrow$
$\rho^{-1}E\beta^\rho_a(X) \in \overline{M_a}^+$
$\Rightarrow$
$\alpha_{-a}\rho^{-1}E\beta^\rho_a(X) \in \overline{M}^+$
$\Rightarrow$
$\rho \alpha_{-a}\rho^{-1}E\beta^\rho_a(X) \in \overline{\rho(M)}^+.$

$E_a$ is an operator valued weight.
We check that for every $X, Y \in M$
\begin{align*}
E_a (\rho(Y)^* X \rho(Y)) & =
\rho \alpha_{-a}\rho^{-1} E ( \rho(\alpha_a(Y^*)) \beta^\rho_a (X)
\rho(\alpha_a(Y)) \\
&= \rho \alpha_{-a} ( \alpha_a(Y^*) \rho^{-1} (E (\beta^\rho_a (X)) )
\alpha_a(Y) )\\
&= \rho  ( Y^* \alpha_{-a} \rho^{-1} (E (\beta^\rho_a (X)) ) Y ) \\
&=\rho(Y^*) \rho( \alpha_{-a} (\rho^{-1} (E (\beta^\rho_a (X)) ))) \rho(Y)\\
&= \rho(Y^*) E_a(X) \rho(Y).
\end{align*}
$E_a$ is normal, thus also faithful and semifinite
(see \cite{St}, p. 155) since
$\omega \rho^{-1} \circ E_a = \omega \rho^{-1} E \beta^\rho_a =
\psi_\rho \beta^\rho_a,$ which is faithful
(since $\psi_\rho$ is)
and semifinite by hypothesis.

In the irreducible case it follows by uniqueness
that $E = \lambda_a E_a,$ $\lambda_a \in {\Reali}_+$,
see e.g. \cite{St}, p.174 Corollary 12.13 (see also \cite{EN}, Prop. 11.1),
thus $\psi_\rho
= \lambda_a \psi_a$
and
$\lambda_a=1$
by Sect. 2, 5.
$\hfill\Box$\medskip

\noindent{\bf 5.10 Corollary.} {\sl
Let $\rho$ be a covariant irreducible morphism of $\A$
localized in $(b,+\infty),$ $b>0.$
Then the following are equivalent:
\begin{description}
\item{1)} $\psi_{a}$ is semifinite for some ($\Leftrightarrow$ for all) $a<b$;
\item{2)} $\rho$ has positive energy;
\item{3)}
$\psi_\rho = \psi_{a}$ for some ($\Leftrightarrow$ for all)
$a>0.$
\end{description}
}\medskip

\noindent{\bf Proof.} It is immediate from
Proposition 2.4, cf. Proposition 5.5, and
Proposition 5.9.
$\hfill\Box$\medskip

Our result entails the following one for sectors on non strongly additive
nets.
\medskip

\noindent{\bf 5.11 Corollary.} {\sl
Let $\A$ be a net on $\Reali$ satisfying the property of Section 3 (but
not necessarily strongly additive). Let $\rho$ be an irreducible covariant
morphism
of $\A$ localized in an bounded interval $I$. If $\rho$ acts identically on
$\A(1,\infty)'\cap\A(0,\infty)$, then $\rho$ has positive energy
iff $\psi_a$ is semifinite.}
\medskip

\noindent{\bf Proof.} We may assume $I=(0,1)$. $\rho$ extends to a morphism
$\tilde\rho$ of the
dual net $\A^d$, covariant with respect to the same representation of
the translation-dilation group. Since
$\A^d(a,0)=\A(0,\infty)'\cap\A(a,\infty)$, $a<0$,
and $\A^d(1,b))=\A(b,\infty)'\cap\A(1,\infty)$, $b>1$,
it follows that $\tilde\rho$ is still localized in $(0,1)$. Now $\A^d$ is
strongly additive \cite{GLW}, thus $\tilde\rho$, hence $\rho$, has positive
energy
by Prop. 5.9.
$\hfill\Box$

\medskip

If $(\rho,V_\rho)$ and $(\sigma,V_\sigma)$ are two covariant morphisms,
then $\rho\sigma$ is covariant with respect to the representation of $P$
given by
$$V_{\rho \sigma}(g):=\rho(z_{\sigma}(g))V_{\rho}(g), g \in P \eqno (5.2)$$
as shown by the following computation ($X \in M, g \in P$):
\begin{align*}
V_{\rho \sigma}(g) (\rho \circ \sigma (X)) V_{\rho \sigma}(g)^*
&= \rho \circ \sigma (V(g)XV(g)^*)
= \rho(V_{\sigma}(g) \sigma(X) V_{\sigma}(g)^*)\\
&= \rho(z_{\sigma}(g)V(g)\sigma (X)V(g)^*z_{\sigma}(g)^*)\\
&= \rho(z_{\sigma}(g))\rho(V(g)\sigma(X)V(g)^*)\rho(z_{\sigma}(g)^*)\\
&= \rho(z_{\sigma}(g))V_{\rho}(g)\rho\circ\sigma(X)
V_{\rho}(g)^*\rho(z_{\sigma}(g))^*.
\end{align*}

In the following we will show that if two
irreducible sectors $\rho, \sigma$ are covariant with positive energy with
respect to the representations $V_\rho, V_\sigma,$ then $V_{\rho\sigma}$
defined above is a positive energy representation.

\medskip
\noindent{\bf 5.12 Proposition.} {\sl
Let $\A$ be a strongly additive net on $\Reali,$ and
$\rho,$ $\sigma$ be two covariant irreducible morphisms with positive
energy.
Then the representation $V_{\rho\sigma}$ defined by equation (5.2)
has positive energy.
}
\medskip

\noindent{\bf Proof.}
To prove this we show the translational invariance of the unbounded left inverse
$\phi_{\rho \sigma}$
(cf. the last subsection), namely
$\alpha^{-1}_a \phi_{\rho \sigma} \beta^{\rho \sigma}_a
=\phi_{\rho \sigma}, \ a >0.$
As preliminary result we shall now prove that:
$$\phi_{\rho\sigma} = \phi_{\sigma}\phi_{\rho} \eqno (5.3)$$
Let $\psi$ be a s.n.f.
weight on the von Neumann algebra $M$ and $\rho$ be an
isomorphism of $M$ onto its subalgebra $\rho(M).$
Then we have
$$\sigma^{\psi\rho^{-1}}_t (Y) = \rho \circ \sigma^{\psi}_t \circ \rho^{-1}(Y),
\ Y \in \rho(M)$$ checking the invariance and the KMS property.
{}From this, given two s.n.f. weights $\psi_i, \ i=1,2,$ on $M$,
we obtain by direct computation
$$\rho ((D\psi_1 : D\psi_2)_t) = (D\psi_1 \rho^{-1} : D\psi_2 \rho^{-1})_t .
\eqno (5.4)$$
Now we have
$$z_{\rho\sigma}(t) = (D\psi_\sigma \rho^{-1}E_\rho : D\omega)_t \eqno (5.5)$$
as shown by the following computation:
\begin{align*}
z_{\rho\sigma}(t) = \rho(z_\sigma(t))z_\rho(t)
&=\rho ((D\psi_\sigma : D\omega)_t) (D\psi_\rho : D\omega)_t \\
\mbox{(by equation (5.4))} \quad\quad &=(D\psi_\sigma \rho^{-1} : D\omega
\rho^{-1})_t
(D\psi_\rho : D\omega)_t \\
\mbox{(by \cite{St} Theorem 11.9)} \quad\quad
&= (D\psi_\sigma \rho^{-1}E_\rho : D\omega \rho^{-1}E_\rho)_t
(D\psi_\rho : D\omega)_t \\
\mbox{(since $\psi_\rho = \omega \rho^{-1}E_\rho$)}  \quad\quad
&= (D\psi_\sigma \rho^{-1}E_\rho : D\psi_\rho)_t
(D\psi_\rho : D\omega)_t \\
\mbox{(by \cite{St} Corollary 3.5)}  \quad\quad
&= (D\psi_\sigma \rho^{-1}E_\rho : D\omega)_t
\end{align*}
{}From Formula (5.5), by a Theorem of Connes (see \cite{St} Corollary 3.6)
we deduce:
$$\psi_\sigma \rho^{-1}E_\rho = \psi_{\rho\sigma},$$
and from a Theorem of Haagerup (see \cite{St} Theorem 11.9):
$$E_{\rho\sigma} = \rho E_\sigma\rho^{-1}E_\rho$$
and finally, applying $(\rho\sigma)^{-1},$ to both sides we get equation (5.3).
The invariance of $\phi_{\rho \sigma}$ is readily obtained using
equation (5.3)
($X \in M^+$):
\begin{align*}
\alpha_{a}^{-1}\circ  &\phi_{\rho \sigma}\circ \beta^{\rho\sigma}_a (X)
= \alpha_{a}^{-1}\circ \phi_{\sigma}\circ\phi_{\rho}\circ
\beta^{\rho\sigma}_a (X)\\
&= \alpha_{a}^{-1}(\phi_{\sigma}\circ\phi_{\rho}\circ \mbox{Ad}(V_{\rho
\sigma}(a))(X))\\
\mbox{by equation (5.2)} \quad
&= \alpha_{a}^{-1}(\phi_{\sigma}\circ\phi_{\rho}\circ
\mbox{Ad}(\rho(z_{\sigma}(a))
V_{\rho}(a))(X))\\
&= \alpha_{a}^{-1}(\phi_{\sigma}(z_{\sigma}(a)\phi_{\rho}(\beta^{\rho}_a(X))
z_{\sigma}(a)^*))\\
&= \alpha_{a}^{-1}(\phi_{\sigma}(z_{\sigma}(a) (\alpha_a\circ \phi_{\rho}(X))
z_{\sigma}(a)^*))\\
&= \alpha_{a}^{-1}(\phi_{\sigma}(\beta^{\sigma}_a \circ\phi_{\rho}(X)))
= \phi_{\sigma}\circ\phi_{\rho}(X).
\end{align*}
The conclusion follows using Prop. 5.5, by the invariance of
$\psi_{\rho \sigma}$.
$\hfill\Box$

\medskip

As already mentioned, the net on $\mathbb R$ we are considering are
obtained by a local
conformal precosheaf by removing a point from the circle \cite{GLW}.
Therefore our results have a version for M\"obius covariant sectors.

\medskip

\noindent{\bf 5.13 Theorem.} {\sl
Let $\A$ be a strongly additive local conformal precosheaf on $S^1.$
The class of M\"obius covariant (resp. traslation covariant, with respect
to a given $\infty$ point)
sectors with positive energy is stable under composition,
conjugation and direct integral decomposition.}
\medskip

\noindent{\bf Proof.}
The stability  of the covariance with positive energy under
direct integral decomposition is shown explicitly in
Lemma 5.14 of the next subsection.

Let thus assume that $\rho, \sigma$ are covariant sectors with positive
energy. Let
$\rho=\int^{\oplus}\rho_{\lambda}\mbox{d}\mu(\lambda)$ and
$\sigma = \int^{\oplus}\sigma_{\lambda'}\mbox{d}\nu(\lambda')$
be two direct integral decompositions into
irreducible sectors.
By the previous statement,
the irreducible components of $\rho$ (resp. $\sigma$) are
$\mu$ (resp. $\nu$) almost everywhere covariant with positive energy.
Then:
$$\rho\sigma=\int^{\oplus}\rho_{\lambda}\sigma_{\lambda'}
\mbox{d}(\mu\times\nu)(\lambda,\lambda')$$
therefore by Proposition 5.12 $\rho_{\lambda}\sigma_{\lambda'}$ is
covariant with positive energy almost everywhere, and the same is
true for $\rho\sigma$ by Lemma 5.14.

It remains to show that if $\rho$ is covariant with positive energy,
the same is true for its conjugate $\bar\rho= j\circ\rho\circ j$,
where $j=\mbox{Ad}J$ and $J$ is the modular conjugation of $(M,\Omega)$.
But this immediately follows by setting
$V_{\bar\rho}(g)=JV_{\rho}(rgr)J$, where $r$ is the change of sign
on $\mathbb R$, see \cite{GuLoa}.
$\hfill\Box$

\subsection{An example of sector with infinite dimension and
negative energy levels}

We show now that there exist translation-dilation covariant sectors
whose energy spectrum is the real line. Our example, concerning a
non strongly additive net, will be a reducible sector, but enlightens
nevertheless the structure of the involved objects and the limit of the
arguments.

\medskip
\noindent{\bf 5.14 Lemma.} {\sl
Let $\rho = \int^{\oplus}\rho_{\lambda}\mbox{d}\mu(\lambda)$ be a
(non unique)
direct integral decomposition
of a sector $\rho$ of
a net of von Neumann algebras $\A$ on $\Reali$  as in section 3 (resp. on
$S^1$).
Then $\rho$ is translation (resp. M\"obius) covariant with
positive energy iff $\rho_{\lambda}$ is translation (resp. M\"obius)
covariant with positive energy $\mu$--almost everywhere.}
\medskip

\noindent{\bf Proof.}
Clearly if
$\rho_{\lambda}$ is translation (resp. M\"obius) covariant with positive
energy for $\mu$-almost
all $\lambda$, the same is true for $\rho$. Conversely if $\rho$ is
translation
covariant with positive energy, there exists by Borchers theorem
\cite{Bo0}
a unitary one-parameter group
$T_{\rho}\in \rho(\A)''$, with positive generator,
implementing the translations $\rho\circ
\mbox{Ad}T(\cdot)$.
Therefore $T_{\rho}$ has a decomposition
$T_{\rho}=\int^{\oplus}T_{\rho}^{(\lambda)}\mbox{d}\mu(\lambda)$,
where $T_{\rho}^{(\lambda)}$ has positive generator for almost all
$\lambda$, and implements the translations on $\rho_{\lambda}$.
If moreover $\rho$ is  covariant with respect to (the universal
covering of the) M\"obius group,
we may repeat the argument with the translations associated with
different intervals of $S^1$ and find implementations in the weak
closure of $\rho_{\lambda}(\A)$. In this way we construct, for almost
all $\lambda$, a unitary
representation $V_{\lambda}$ of the universal covering of the
M\"obius group up to a cocycle in the
center of $\rho_{\lambda}(\A)'$. Since the cohomology of the universal
covering of the M\"obius group is trivial, we may replace
$V_{\lambda}(g)$ with $z(g)V_{\lambda}$, where $z(g)$ in the center of
$\rho(\A)'$, and get a true unitary representation providing the
M\"obius covariance with positive energy of $\rho_{\lambda}$.
$\hfill\Box$
\medskip

\noindent{\bf 5.15 Proposition.} {\sl
Let $\A$ be an irreducible net of von Neumann algebras on $\Reali,$ covariant
under translations and dilations as in Sect. 3.
Let $\gamma$ be a morphism of $\A$ localized in the interval
$I \subset \Reali,$ and assume that
$\gamma_a:=\alpha_{-a}\gamma\alpha_a$
is disjoint from $\gamma$ for every $a \neq 0$
(thus $\gamma$ is not translationally covariant).
Then $\rho:= \int^{\oplus} \gamma_a da$
is a translationally covariant (reducible) endomorphism
with infinite statistics whose energy spectrum is $\mathbb R$.
}
\medskip

\noindent{\bf Proof.} $\rho$ acts on
vectors $\xi$ in
the Hilbert space $L^2({\Reali},\H)$ (separable if $\H$ is separable)
via $(\rho(X)\xi)(a) = \gamma_{-a} (X) \xi(a),$ $X \in \A,$ $a \in \Reali,$
and covariance under translations is implemented by the unitary one--parameter
group $(T(b)\xi)(a)=\xi(b-a).$
However $\rho$ has not positive energy since
otherwise we would infer from Lemma 5.14 that
(for almost every $a \in \Reali$) $\gamma_a$ is covariant and this
is not possible by hypothesis.
$\hfill\Box$
\medskip

To give an explicit example, we recall  that
although
the free scalar massless field
$\varphi = \varphi(t,x)$ in two dimensions does not exist,
its derivative $j=\partial_0 \varphi - \partial_1 \varphi$
makes sense and depends on $t-x$, thus defines a one dimensional
current $j$.
Every $f \in {\cal S}(\Reali)$
such that $\int_{\Reali} f(t) dt = 0$ determines
a unitary operator $W(f)=e^{ij(f)}$ and the Weyl relations
$W(f + g) = e^{i \int fg' dt}W(f)W(g)$ are satisfied.
Let $\A$ be the (strongly additive) net on $\Reali$
defined by
$\A(I):=\{ W(f) \ | \ f \in {\cal S},\ \int_{\Reali} fdt=0,
\ \mbox{supp}(f) \subset I\}.$
As is known \cite{BMT}, this net has localized automorphisms $\alpha_q,$
with $q \in {\cal S}$ real valued with compact support,
given by $\gamma_q(W(f))=e^{2i \int q f} W(f)$
\cite{BMT}.
If $q = Q'$ with $Q \in {\cal S},$ that is $\int q(t) dt = 0,$
then $\gamma_q$ is inner, indeed $\gamma_q = \mbox{Ad}(W(Q)).$
The equivalence class of these automorphisms are labeled by
the real numbers $q_0 := \int q(t)dt.$

As noticed in \cite{GLW},
it is possible to generalize this construction in order to
obtain non-covariant automorphisms, thus sectors
with the properties needed in Prop. 5.15.
For the sake of completeness we briefly recall this construction.

Let $\B \subset \A$ the net generated by the derivative of $j$:
$\B(I):=\{ W(f) \ | \ f = F ', \ F \in {\cal S},
\ \int F dt = 0, \ \mbox{supp}(F) \subset I \}.$
Then $\B$ is a proper subnet of $\A$ which is not strongly additive,
but $\A(I)=\B(I)$ if $I$ is a half-line. If $q$ is a smooth function
on $\mathbb R$ such that $q'$ has compact support, then $\gamma_q$
makes sense as automorphism of $\B$ and its equivalence class is
labeled by the charges $\int q'(t)dt$ and $\int tq'(t)dt$. As a
consequence, if $q(+\infty)\neq q(-\infty)$, then $\gamma_q$ is a
transportable localized automorphism of $\B$ such that
$\alpha_{-a}\gamma_q \alpha_a$ is disjoint from $\gamma_q$ for each
non trivial translation $\alpha_a.$

\subsection{Construction of the M\"obius covariant unbounded left inverse}

Let $\A$ be a local conformal precosheaf on $S^1$, namely a map
$$I\to\A(I)$$ that associates to each (proper) interval of $S^{1}$ a
von Neumann algebra, satisfying isotony, locality, M\"obius
covariance
with positive energy, uniqueness of the vacuum, see e.g.
\cite{GuLoc}.
A morphism $\rho$ localized in the interval $I_0$ is a map
$$ I \to \rho_I $$
which associates to every
interval $I$ of $S^1$ a normal representation of
$\A(I)$ on $\H$ such that
$$ \rho_{\tilde{I}}|_{\A(I)} = \rho_I, \ I \subset \tilde{I} $$
and
$$\rho_{I'_0} = \mbox{id}.$$
By Haag duality $\rho_I \in \mbox{End}(\A(I))$ if $I \supset I_0.$

We now assume that $\A$ is strongly additive. Then $\rho$ is automatically
covariant with positive energy if $d(\rho) < \infty.$
In the following $d(\rho)$ is infinite and we assume M\"obius
covariance with positive energy.

By an {\it unbounded left inverse}
$\phi$ of $\rho$ we shall mean a map
$I \to \phi_I$ that associates with any interval $I \supset I_0$
a map $\phi_I: \ \A(I)^+ \to \overline{\A(I)}^+$ such that
$$ \phi_I ( \rho_I(u) X \rho_I(u^*) ) = u \phi_I(X) u^*, \ u, X \in \A(I)$$
i.e. $\rho_I \phi_I$ is a $\rho(\A(I))$--valued weight on $\A(I),$
and
$$
\phi_{\tilde I}|_{\A(I)} = \phi_I
$$
if $I\subset \tilde I$ are intervals containing $I_0$.

\medskip
\noindent{\bf 5.16 Proposition.} {\sl
Let $\A$ be strongly additive and $\rho$
irreducible,
covariant with positive energy and
localized in $I_0.$
There exists an unbounded left inverse $\phi$ of $\rho,$
covariant with respect to the
M\"obius group, namely
such that
$$ \alpha_g^{-1} \circ \phi_{gI} \circ \beta^\rho_g \ = \ \phi_{I}$$
whenever $I \cap gI \supset I_0, \ g \in P.$
}\medskip

\noindent{\bf Proof.}
Let $\psi_I = \psi_{\rho,I}$ be the s.n.f. weight on $\A(I),$
$I \supset I_0,$
defined by
$(D\psi_I : D\omega_I)_t = z_\rho(-2 \pi t), \ t \in \Reali.$
Then
$$\psi_{\tilde{I}}|_{\A(I)} = \psi_I, \ I \subset \tilde{I} \eqno (5.6).$$
To check this notice that (5.6) is true if $I$ and $\tilde{I}$
have a common boundary point, as by
cutting the circle we may assume
$I=(1,\infty),$ $\tilde{I} = (0,\infty)$ and we may apply
the results so far obtained
for the real line.
Thus (5.6) is true in general as we may check it in two steps
by considering an intermediate interval
$I \subset I_1 \subset \tilde{I}$ such that
both $I \subset I_1$ and $I_1 \subset \tilde{I}$
have a common boundary point.\par
Concerning the invariance, let $I_1, I_2 \subset \tilde{I}$
be intervals
containing $I_0$
and $g \in \mbox{\tt PSL}(2,\Reali)$ with $g I_1 = I_2.$
If $X_2 \in \A(I_2)$ then $X_2 = \beta^\rho_g(X_1)$ with $X_1 \in \A(I_1).$

In fact $z_\rho(g) \in \A(g I_1)$
since $I_0$ and $g I_0$ are both contained in $g I_1.$
Then
\begin{align*}
\psi_{gI_1}( \beta^\rho_g(X_1) ) & = \psi_{I_2} (X_2) = \psi_{\tilde{I}}(X_2) \\
& =  \psi_{\tilde{I}} ( \beta^\rho_g(X_1) ) = \psi_{\tilde{I}} (X_1)\\
& = \psi_{I_1}(X_1)
\end{align*}

\hfill (5.7)

\noindent
provided $g$ is a translation or dilation of $\tilde{I}.$

It follows that
(5.7) holds for all $g$ such that
$g I, I \supset I_0 $
by a simple repetition of the arguments.
The reason is as follows:
given two such intervals $I, gI$ we can always find sequences
$I_1 = I, \ldots, I_n = gI$ and $\{ \tilde{I}_i \}_{1}^{n-1}$
(actually $n=3$)
such that
$I_i \cap I_{i+1} \supset I_0,$ $I_i \cup I_{i+1} \subset \tilde{I}_i,$
$I_{i+1} = g_i I_i$ with $g_i$ a translation of $\tilde{I}_i .$
Namely, if one of the two intervals $I_1=I, I_3=gI$ is contained in the other,
let us
take as $I_2$ an intermediate interval containing the smallest of the two and
one of the endpoints of the biggest.
Otherwise take as $I_2$ the connected component of $I_0$ in $I \cap gI.$
In both of the situations it is easy to see that $I_2$ is an interval included
in one of the two intervals $I_1, I_3,$
containing the other one
and having the two endpoints in common with $I_1$
and $I_3$ respectively.
It is now possible to apply the previous results in two successive steps to
the pairs $I_1, I_2$ and $I_2, I_3$ taking as
$\tilde{I}_1$ and $\tilde{I}_2$
the biggest intervals and considering the unique translation
$g_1$ of $\tilde{I}_1$ such that $g_1 I_1 = I_2$
(resp. $g_2$ of $\tilde{I}_2$ such that $g_2 I_2 = I_3$)
with fixed point the common extreme of $I_1$ and $I_2$ (resp. $I_2$ and $I_3$).
Then $g=g_1 g_2 h$ where $h$ is a dilation of $I$ and we have (by (5.7)):
$\psi_{g_1g_2hI}(\beta^\rho_{g_1}\beta^\rho_{g_2}\beta^\rho_{h}(X))=
\psi_I(X),$ $X \in \A(I).$
Let $E_I : \A(I)^+ \to \overline{\rho_I(\A(I))}^+$ the Haagerup's
operator valued weight, then
$$ \phi_I := \rho^{-1}_I E_I$$
is the desired unbounded left inverse.
$\hfill\Box$\medskip

\medskip

\noindent{\bf 5.17 Corollary.} {\sl
Let $\A$ be strongly additive and $\rho$
irreducible,
covariant and
localized in $I_0.$
Then $\rho$ has positive energy if and only if there exists an unbounded
left inverse
$\phi$ of $\rho,$
covariant with respect to the
M\"obius group.
}\medskip

\noindent{\bf Proof.}
The existence of the M\"obius covariant unbounded left inverse has
been shown in the preceding Proposition 5.16.
The reverse implication follows from the invariance of
$\psi$ as in Corollary 5.10.
$\hfill\Box$\medskip

\medskip
\appendix
\section{Appendix. Araki relative modular operators and Connes spatial
derivatives}\label{Appendix A}

We use the notation in \cite{St}, Ch. I, II.
In particular given a von Neumann algebra $M \subset B(\H),$ a
semifinite normal faithful (s.n.f.) weight $\psi$
on $M$,
$\I_\psi$ will denote the dense left ideal
$\{X \in M \ | \ \psi(X^*X) < \infty \}$,
$\H_\psi$ the GNS Hilbert space of $\psi$. $\I_\psi$
is embedded as a dense linear subspace of $\H_\psi$, denoted by
$X\to (X)_{\psi}$, and the GNS representation
$\pi_{\psi}$ is given by $\pi_{\psi}(X)(Y)_{\psi}
=(XY)_{\psi}$, $X\in M,\ Y\in\I_\psi$.
$\D(\H,\psi)$ will denote the dense linear subspace of all
$\zeta \in \H$ that are $\psi$-bounded,  namely
the linear operator $R^\psi_\zeta: \H_\psi \to \H$  defined by
$(X)_\psi \to X\zeta, \ X \in \I_\psi$ is bounded.
If $\chi'$ is a s.n.f. weight on $M'$,
$\Delta(\psi /\chi')$ denotes the
Connes spatial derivative.\footnote{If $A$ is a
positive linear operator,
we set $\|A\zeta\|=+\infty$ for all vectors $\zeta\notin
\D(A)$.}
\medskip

\noindent{\bf A.1 Lemma.} {\sl
Let $M \subset B(\H)$ be a von Neumann algebra,
$\omega'=(\Omega, \cdot \Omega)$ a
vector state on $M'$ given by a cyclic and separating vector $\Omega,$
$\psi$ a normal faithful semifinite weight on $M$.
Then
$\|\Delta(\psi /\omega')^{\frac{1}{2}}X\Omega \|^2 =
 \psi(XX^*), \ X \in M.$
}\medskip

\noindent{\bf Proof.}
As $\omega'$ is the state on $M'$ given by a vector $\Omega$, one
checks immediately that
$M\Omega = \D(\H,\omega')$
and $R^{\omega'}_{X\Omega} = X.$
The Lemma is thus obtained by taking $\zeta = X \Omega$ in the
formula $\|\Delta(\psi /\omega')^{\frac{1}{2}} \zeta \|^2 =
\psi(R^{\omega'}_\zeta (R^{\omega'}_\zeta)^*)$ for $\zeta \in \D(\H,\omega')$
that defines the spatial derivative \cite{St} p. 95.
$\hfill\Box$\medskip

By polarization we also deduce that
$$(\Delta(\psi /\omega')^{\frac{1}{2}}Y\Omega,
\Delta(\psi/\omega')^{\frac{1}{2}}X\Omega )
= \psi(XY^*), \ X, Y \in \I^*_\psi.$$
We also need two facts.
Recall (see \cite{St} Chapters I, II) that if $\psi$
and $\omega$ are s.n.f. weight on $M$ then the
Tomita operator
$S_{\psi,\omega}: \H_\omega \to \H_\psi$
is the closure of the densely defined anti-linear
operator
$(X)_{\omega}\to (X^*)_\psi,$
$X \in \I_{\omega}\cap\I^*_\psi.$
The polar decomposition
$S_{\psi, \omega} = J_{\psi, \omega} \Delta_{\psi, \omega}^{\frac{1}{2}},$
$\Delta_{\psi ,\omega}^{\frac{1}{2}} =
(S^*_{\psi,\omega} S_{\psi,\omega})^{\frac{1}{2}},$
defines the relative modular conjugation
and the relative modular
operator.

Let  $V_{\psi, \omega}: \H_\omega \to \H_\psi$ be the
uniquely determined unitary operator satisfying:
$\pi_\psi (X) = V_{\psi, \omega}\pi_\omega (X) V_{\psi, \omega}^*,$
preserving the natural cones,
see \cite{St}, 3.16. We shall identify $\H_{\omega}$
and $\H_\psi$ via $V_{\psi,\omega}$, cf. \cite{St}, 3.17,
therefore $V_{\psi, \omega}=1$ and $J_{\psi, \omega}
= J_{\omega , \omega}$, that we simply denote by $J$
in the following. We will assume $M$ to act standardly
on $\H_{\psi}=\H_{\omega}$, thus we suppress the symbols $\pi_{\psi}$
and $\pi_{\omega}$, and
shall consider the s.n.f.
weight $\omega'= \omega(J\cdot J)$ on $M' .$

\medskip

\noindent{\bf A.2 Lemma.}
{\sl
Let $M,$ $\psi,$ $\omega'$ as above; then
$$\Delta(\psi /\omega')^{\frac{1}{2}}
= {\Delta^{\frac{1}{2}}_{\psi,\omega}}.$$}
\par
\noindent{\bf Proof.}
Let us first assume that $M$ is a factor. Then
$\Delta(\psi /\omega')^{it}
{\Delta^{-it}_{\psi,\omega}} \in M \cap M'=\mathbb C,$ $t \in \Reali$
because both $\Delta(\psi /\omega')^{it}$ and
${\Delta^{it}_{\psi,\omega}}$ implement the same modular groups
of $\omega'$ on $M'$ and of $\psi$ on $M$.
It follows that
$\Delta(\psi /\omega') =
\lambda\Delta_{\psi,\omega}$ for some $\lambda>0.$
Now, for every $X \in \I_\omega \cap \I^*_\psi,$
we have (cf. the proof of Lemma A.1)
\begin{align*}
\| \Delta^{\frac{1}{2}}_{\psi,\omega} (X)_\omega \|^2
&= \| J_{\psi,\omega} \Delta^{\frac{1}{2}}_{\psi,\omega} (X)_\omega \|^2
= \|(X^*)_\psi \|^2 \\
&= \psi(XX^*)
= \psi(R^{\omega'}_{(X)_\omega} {R^{\omega'}_{(X)_\omega}}^*) \\
&=\|\Delta(\psi /\omega')^{\frac{1}{2}} (X)_\omega \|^2 .
\end{align*}

\hfill (A.1)

In fact
from the relation
$$JYJ (X)_\omega = X J (Y)_\omega, \ Y \in \I_\omega, \eqno (A.2)$$
see \cite{St}, p. 26,
it follows that $(X)_\omega \in \D(\H_\omega,\omega')$
since
\begin{align*}
(JYJ (X)_\omega, JYJ (X)_\omega)_\omega
&= (X J (Y)_\omega,X J (Y)_\omega)_\omega
= (J (Y)_\omega, X^*X J (Y)_\omega)_\omega \\
& \leq \| X \|^2 ((Y)_\omega,(Y)_\omega)_\omega
= \| X \|^2 \omega(Y^*Y) \\
&=\| X \|^2 \omega' ((JYJ)^*(JYJ));
\end{align*}
furthermore
$R^{\omega'}_{(X)_\omega} =X$
as results from:
\begin{align*}
R^{\omega '}_{(X)_\omega}(JYJ)_{\omega '} &= JYJ(X)_{\omega} \\
(\mbox{by equation (A.2)}) \quad &= XJ(Y)_{\omega}\\
&=X(JYJ)_{\omega '}
\end{align*}
where, as before, $Y \in \I_{\omega}$
and we have identified $(JYJ)_{\omega'}$ and $J (Y)_\omega$.
Thus by equation (A.1), we obtain
$\lambda=1.$
The conclusion follows by direct integral
decomposition.
$\hfill\Box$
\medskip

If $A$ is a positive operator on $\H$ affiliated to $M$ we
define
$$\psi(A):= \mbox{sup}_n \psi(A E_A ([0,n)) ),$$
where $E_A$ is the projection--valued spectral measure of $A$.

\medskip

\noindent{\bf A.3 Lemma.}
{\sl
Let $M \subset B(\H),$ $\psi,$ $\omega',$ $\Omega$ as above; then
$$
\D (\Delta(\psi /\omega')^{\frac{1}{2}}) =
\{ T \Omega \ | \ T \ \mbox{closed}, \ T \ \eta \ M,
\ \Omega \in \D (T),
\ \psi(TT^*) < \infty \}.
$$}
\par
\noindent{\bf Proof.}
$\subset \ :$
given $\zeta \in \D (\Delta(\psi /\omega')^{\frac{1}{2}}),$
consider the densely defined operator
$T^0_\zeta : X' \Omega \to X' \zeta , \ X' \in M'$.
Then  $T^0_\zeta $ is closable since its adjoint is densely defined;
indeed
$J(Y)_\psi \in \D ({T^0_\zeta}^*),$ $Y \in \I_\psi,$
and
${T^0_\zeta}^* J (Y)_\psi =
J Y {\Delta_{\psi,\omega}}^{\frac{1}{2}} \zeta$
since
\begin{align*}
(T^0_\zeta JX\Omega, J&(Y)_\psi)
=(JXJ\zeta, J(Y)_\psi)
=(\zeta, JX^*J J(Y)_\psi)\\
&=(\zeta, J X^*(Y)_\psi)
=(\zeta, J(X^* Y)_\psi)\\
&=(\zeta, JS_{\psi,\omega}Y^* X\Omega)
=(\zeta,\Delta_{\psi,\omega}^{\frac{1}{2}}Y^* X\Omega)\\
&=(\Delta_{\psi,\omega}^{\frac{1}{2}}\zeta,Y^* X\Omega)
=((Y\Delta_{\psi,\omega}^{\frac{1}{2}}\zeta,X\Omega)
=(JX\Omega,JY\Delta_{\psi,\omega}^{\frac{1}{2}}\zeta).
\end{align*}
Furthermore $T^0_\zeta$ commutes with unitaries in $M'$, thus its closure
$T_\zeta$ is affiliated to $M.$
Let $T_\zeta = V H$ be the polar decomposition,
and $e_n : = E_{H} ([0,n))$, so that
$T_n:= T_\zeta e_n \in M,$ and $T_n \Omega \to T_\zeta \Omega =\zeta.$
Finally for each $X \in \I^*_\psi$ we have
$$(T_n \Omega, S^*_{\psi,\omega} JX \Omega)=
(J X \Omega, e_n J {\Delta_{\psi,\omega}}^{\frac{1}{2}} \zeta).$$
This is shown, using the fact that
\begin {align*}
S_{\psi, \omega}^* &(JX \Omega)
=J S_{\psi, \omega}X\Omega\\
&=JS_{\psi, \omega}JX\Omega
=J(X^*)_\psi,
\end{align*}

by the following computation:
\begin{align*}
(T_\zeta e_n \Omega, S^*_{\psi, \omega}JX \Omega)
&=(T_\zeta e_n \Omega, J(X^*)_\psi)\\
&=(T_\zeta e_n \Omega, J(X^*)_\psi)=
(e_n \Omega, T^*_\zeta J(X^*)_\psi)\\
&=(e_n \Omega, J X^* \Delta_{\psi, \omega}^{\frac{1}{2}} \zeta)=
(\Omega, e_n J X^* \Delta_{\psi, \omega}^{\frac{1}{2}} \zeta)\\
&=(J e_n J X^* \Delta_{\psi, \omega}^{\frac{1}{2}} \zeta,
\Omega)
=(J e_n J \Delta_{\psi, \omega}^{\frac{1}{2}}\zeta,
X \Omega)\\
&=(J X \Omega, e_n J \Delta_{\psi, \omega}^{\frac{1}{2}}
\zeta)
=(J X \Omega,  e_n J
\Delta_{\psi, \omega}^{\frac{1}{2}}\zeta),
\end{align*}

therefore $T_n \Omega \in \D (\Delta(\psi /\omega')^{\frac{1}{2}})$
and
$$S_{\psi, \omega}T_\zeta e_n \Omega= e_n J
\Delta_{\psi, \omega}^{\frac{1}{2}}\zeta.
$$
Moreover
$$\psi(T_n T^*_n) =
\| \Delta(\psi /\omega')^{\frac{1}{2}} T_n \Omega \|^2 \leq
\| \Delta(\psi /\omega')^{\frac{1}{2}} \zeta \|^2 .$$
Finally $\psi(T_n T^*_n) \nearrow \psi(T_\zeta T^*_\zeta) < \infty.$
\par
$\supset \ :$
let $T_n := e_n T \in M,$
where $e_n$ is defined using the spectral family of $TT^*;$
then
$T_n T^*_n = e_n TT^* $
is increasing sequence of (bounded) positive operators,
$T_n T^*_n \leq TT^*,$
$\psi(T_n T^*_n) =
\| \Delta(\psi /\omega')^{\frac{1}{2}} T_n \Omega \|^2 \leq
\psi(TT^*) < \infty.$
Therefore $T_n \Omega =e_n T \Omega \to \zeta,$
and
$\Delta(\psi /\omega')^{\frac{1}{2}} T_n \Omega$ is convergent
i.e. it is a Cauchy sequence since
$\Delta(\psi /\omega')^{\frac{1}{2}} (T_n - T_m) \Omega=
\psi((T_n - T_m)(T_n - T_m)^*) = \psi((e_n - e_m)TT^*) =
\psi(e_n TT^*) - \psi(e_m TT^*) \to 0$ when $n, m \to \infty.$
In particular it follows that
$\zeta \in \D (\Delta(\psi /\omega')^{\frac{1}{2}}).$
$\hfill\Box$\medskip

Before concluding this appendix, we recall the notion of the
extended positive part of a von Neumann algebra, needed in Sect. 6.

Given a von Neumann algebra $M$,
its extended positive part $\overline{M}^+$
is defined as the family of all
additive, positively homogeneous and lower semicontinuous
functions $m: M_*^+ \to [0,\infty]$
 \cite{St}, 11.1.
$\overline{M}^+$ has the following characterization  (see \cite{St},
11.3):
any element $m \in \overline{M}^+$ can be uniquely represented as
$m= he+(1-e)\infty$  where $e$ is
a projection in $ M$ and $h$ is a positive self--adjoint operator $h$
affiliated to $eMe$.
Every normal weight on $M$ has a canonical extension
to $\overline{M}^+$ such that $\varphi(m)=m(\varphi), \ \varphi \in M_*^+$
(see \cite{St}, 11.4).

\enddocument